\documentclass[preprint, 12pt]{elsarticle}
\usepackage{graphicx} 
\usepackage{amsfonts,amsmath,amssymb,amsthm}
\usepackage{hyperref}
\usepackage{mathtools}
\usepackage{nicematrix}
\usepackage{float}
\usepackage{caption}
\usepackage{subcaption}
\numberwithin{equation}{section}
\allowdisplaybreaks

\begin{document}
\title{\textbf{Stability of Edelstein Hidden-Variable Fractal Interpolation Functions}}
\author[1]{Aiswarya T\corref{cor1}}%
\ead{aiswaryasidhu8113@gmail.com}
\author[1]{Srijanani Anurag Prasad}
\ead{srijanani@iittp.ac.in}
\cortext[cor1]{Corresponding author}
\date{}
\fntext[fn1]{Department of Mathematics and Statistics, IIT Tirupati, Yerpedu P.O., India, 517619}
\newtheorem{proposition}{Proposition}
\newtheorem{theorem}{Theorem}
\newtheorem{definition}{Definition}
\newtheorem{example}{Example}
\newtheorem{remark}{Remark}
\newtheorem{corollary}{Corollary}
\newtheorem{lemma}{Lemma}

\begin{abstract}
 Quantitative stability of Edelstein hidden-variable fractal interpolation functions (EHVFIFs) and the upper box-counting dimension of their graphs is investigated under perturbations of  interpolation nodes, ordinate data, hidden-variable data, and vertical scaling functions.
\end{abstract}
\begin{keyword}
     Edelstein Contraction \sep Fractal Interpolation Function \sep Stability \sep Box-counting Dimension
    \hspace{2mm}
\end{keyword}
\maketitle

\section{Introduction}
Fractal interpolation functions (FIFs), introduced by Barnsley~\cite{barnsley1986fractal}, provide a powerful method for constructing continuous functions whose graphs exhibit fractal characteristics. Since their introduction, fractal interpolation functions have been extensively studied due to their capability to model irregular phenomena and complex data sets in various areas such as approximation theory, signal processing, and image reconstruction~\cite{raubitzek2021fractal, manousopoulos2011parameter, banerjee2020multifractal, chen2011reconstruction}. Over the years, several extensions of classical fractal interpolation functions have been proposed to enhance their flexibility and approximation capabilities. One significant generalisation is the hidden-variable fractal interpolation function (HVFIF)~\cite{barnsley1989hidden}, which introduces an additional hidden variable into the interpolation scheme to allow greater control over the geometry and smoothness of the resulting interpolant, thereby leading to richer classes of fractal functions. 

The construction of fractal interpolation functions typically relies on fixed-point theorems associated with contraction mappings. In this context, the contraction principle proposed by Edelstein in~\cite{Edelstein1962} provides a useful generalisation of the classical Banach contraction principle. Edelstein contractions relax the strict Lipschitz contractivity condition while still ensuring the existence of a unique fixed point under suitable assumptions. This generalised framework enables the construction of a broader class of fractal interpolation functions, referred to as Edelstein fractal interpolation functions. More recently, the combination of hidden-variable techniques with Edelstein contractions has led to the development of Edelstein hidden-variable fractal interpolation function(EHVFIF) in~\cite{aiswarya2026}. These functions integrate the advantages of both frameworks: the geometric flexibility provided by hidden-variables and the generalised contractive structure arising from Edelstein mappings. Such constructions yield richer classes of fractal interpolants and open new avenues for studying their analytical and approximation properties.

\section{Preliminaries}\label{prelims}

Let $V=\{(t_j,u_j)\in I\times\mathbb{R}: j=0,1,\ldots,N\},$ be a given data set   in which $\{t_j\}_{j=0}^N$ is strictly increasing and  $V^*=\{(t_j,u_j,v_j)\in I\times\mathbb{R}^2: j=0,1,\ldots,N\}$ denote the corresponding generalised data set. Let $I=[t_0,t_N]$, $I_j=[t_{j-1},t_j]$ for $j=1,\ldots,N$, and  $L_j:I\to I_j$ be contractive homeomorphisms given by 
\begin{align}\label{Ljdefn}
L_j(t)= \frac{t_j-t_{j-1}}{t_N-t_0}  (t - t_0) + t_{j-1}.
\end{align}
Let $F_j:I\times\mathbb{R}^2\to\mathbb{R}^2$ be functions given by
\begin{align}\label{F_jdefn}
F_j(t,u,v)=D_j(t)S_j(u,v)+Q_j(t),
\end{align}
where
\renewcommand{\arraystretch}{1.3}
\[ 
D_j(t)= \begin{bmatrix} \alpha_j(t)& \beta_j(t)\\ \gamma_j(t)& \delta_j(t) \end{bmatrix},\quad
S_j(u,v)= \begin{bmatrix} s_j(u)\\ r_j(v) \end{bmatrix},\quad
Q_j(t)= \begin{bmatrix} p_j(t)\\ q_j(t) \end{bmatrix}, \] 
in which $\alpha_j,\beta_j,\gamma_j,\delta_j, p_j $ and $q_j$  are Lipschitz functions such that \\ $\|\alpha_j\|_\infty +\|\gamma_j\|_\infty\leq 1, \quad \|\beta_j\|_\infty+\|\delta_j\|_\infty\leq 1,\quad F_j(t_0,u_0,v_0)=(u_{j-1},v_{j-1})$ and $F_j(t_N,u_N,v_N)=(u_j,v_j).$ Also, $s_j$ and $r_j$ are real valued Edelstein contractions on $\mathbb{R}$.   Thus, $F_j$ are Edelstein contractions with respect to the Manhattan metric in the second and third variables. The iterated function system (IFS) is defined as $\mathcal{I}=\{I\times K, \omega_j: j=1,\ldots,N\},$ where $K\subset\mathbb{R}^2$ is a compact set such that $F_j(I\times K)\subset K$ and $\omega_j(t,u,v)=(L_j(t),F_j(t,u,v)).$ Suppose $\mathcal{C}(I,\mathbb{R}^2)$ denote the space of $\mathbb{R}^2$-valued continuous functions on $I$ equipped with Manhattan metric,  $\mathcal{C}_e(I,\mathbb{R}^2)=\{h\in \mathcal{C}(I,\mathbb{R}^2):h(t_0)=(u_0,v_0),\,h(t_N)=(u_N,v_N)\}$ and $\mathcal{C}_d(I,\mathbb{R}^2)=\{h\in \mathcal{C}_e(I,\mathbb{R}^2):h(t_j)=(u_j,v_j),\ j=1,\ldots,N-1\}.$  

\begin{proposition}~\cite{aiswarya2026}\label{RBoperatorthm}
The Read–Bajraktarevi\'{c} operator $\mathcal{R}:\mathcal{C}(I,\mathbb{R}^2)\to \mathcal{C}(I,\mathbb{R}^2)$ defined by $\mathcal{R}h(t)=F_j(L_j^{-1}(t),h(L_j^{-1}(t))), \ \mbox{for}\ t\in I_j$ maps an element of $\mathcal{C}_e(I,\mathbb{R}^2)$ to $\mathcal{C}_d(I,\mathbb{R}^2)$ and admits a unique fixed point $f\in \mathcal{C}_d(I,\mathbb{R}^2)$. Moreover, the mappings $\omega_j$ are Edelstein contractions with respect to a metric equivalent to the Euclidean metric, and the graph of $f$ is the attractor of the IFS $\mathcal{I}$.
\end{proposition}

\begin{definition}~\cite{aiswarya2026}
The first component $f_1$ of the vector-valued fixed point $f$ is called the \emph{Edelstein hidden-variable Fractal Interpolation Function (EHVFIF)} associated with the data set $\{(t_j,u_j):j=0,1,\ldots, N\}$.
\end{definition}
\section{Stability of the EHVFIF}\label{stability}

\subsection{Perturbation in the data set}

Let $I=[0,1]=[t_0,t_N]$. Let $f_1$ be the EHVFIF constructed corresponding to the data set $V=\{(t_j,u_j):j=0,1,\ldots,N\}$ and the generalised data set $V^{*}=\{(t_j,u_j,v_j):j=0,1,\ldots,N\}$ with respect to the iterated function system
\begin{align}
    \mathcal{I}=\{I\times K,\ \omega_j : j=1,2,\ldots,N\},\label{IFSV*}
\end{align}
where $\omega_j(t,u,v)=(L_j(t),F_j(t,u,v))$ in which the mappings $L_j$ and $F_j$ are defined as in Eq.~\eqref{Ljdefn} and Eq.~\eqref{F_jdefn}, respectively.

Let $V_1=\{(t_j^{'},u_j):j=0,1,\ldots,N\}$ be another data set on $I$ such that $0=t_0^{'}<t_1^{'}<\ldots<t_N^{'}=1$ which differs from $V$ in the first coordinate and $V_1^{*}=\{(t_j^{'},u_j,v_j):j=0,1,\ldots,N\}$ be the generalised data set corresponding to $V_1$. 
 \begin{theorem}\label{f'=fUinv}
Let $V^*=\{(t_j, u_j, v_j):j=0,1,\ldots,N\}$  and $V_1^*=\{(t_j^{'}, u_j, v_j):j=0,1,\ldots,N\}$  be two generalised data sets, where $t_0=t_0'=0$ and $t_N=t_N'=1$. Then, $f$ is the vector-valued FIF corresponding to the IFS $\mathcal{I}$ if and only if $f\circ U^{-1}$ is the vector-valued FIF corresponding to the IFS $\mathcal{I^{'}}$.
\end{theorem}
 \begin{theorem}\label{tperturbthm}
Let $f_1$ and $f_1^{t}$ be the EHVFIFs corresponding to the generalised data sets $V^{*}$ and $V_1^{*}$, respectively. Then, 
\begin{align}
    \|f_1-f_1^{t}\|_\infty<K_1\max\limits_{1\leq k\leq N}|t_k-t_k^{'}|^\alpha.
\end{align}
\end{theorem}

\begin{theorem}\label{uperturbthm}
Let $f_1$ and $f_1^{u}$ be the EHVFIFs corresponding to the generalised data sets $V^{*}=\{(t_j,u_j,v_j):j=0,1,\ldots,N\}$ and $V_2^{*}=\{(t_j,u_j',v_j):j=0,1,\ldots,N\}$, respectively. Suppose $\zeta+\eta<1$. Then,
\begin{equation}\label{xperturb}
    \|f_1-f_1^{u}\|_{\infty}< 3 \Big(\frac{1+\zeta-\eta}{1-\zeta-\eta}\Big)\max_{j=0,1,\ldots,N}\{|u_j-u_j^{'}|\},
\end{equation}
where $\zeta$ and $\eta$ are defined as earlier.
\end{theorem}

\begin{theorem}\label{vperturbthm}
    Let $f_1$ and $f_1^{v}$ be the EHVFIFs corresponding to the generalised data sets $V^{*}$ and $V_3^{*}$, respectively. Let $\zeta$ and $\eta$ be as defined earlier and suppose that $\zeta+\eta<1.$ Then,
\begin{equation}\label{yperturb}
    \|f_1-f_1^{v}\|_{\infty}<\frac{6\zeta}{1-\zeta-\eta}\max_{j=0,1,\ldots, N}\{|v_j-v_j^{'}|\}.
\end{equation}
\end{theorem}
\begin{theorem}
Let $f_1$ and $\hat{f_1}$ be the EHVFIFs corresponding to the data sets $V^{*}=\{(t_j,u_j,v_j):j=0,1,\ldots ,N\}$ and $\tilde{V}^{*}=\{(t_j^{'},u_j^{'},v_j^{'}):j=0,1,\ldots ,N\}$, respectively. Then,
\begin{equation}\label{allperturb}
\|f_1-\hat{f_1}\|_{\infty}<\rho(V^{*},\tilde{V}^{*}).
\end{equation}
\end{theorem}

\subsection{General case}

\begin{theorem}
Let $f=(f_1,f_2)$ and $f^{u}=(f_1^{u},f_2^{u})$ be the vector-valued FIF corresponding to the generalised data sets $V^{*}=\{(t_j,u_j,v_j):j=0,1,\ldots, N\}$ and $V_2^{*}=\{(t_j,u_j^{'},v_j):j=0,1,\ldots, N\}$, respectively. Then,
\begin{align}
\nonumber\|f_1-f_1^{u}\|_{\infty}&< \frac{1-\eta}{1-\zeta-\eta}\ \max_{j=1,2,\ldots,N}\{l_{p_j}+l_{p_j^{'}}\}+\frac{\zeta}{1-\zeta-\eta}\ \max_{j=1,2,\ldots,N}\{l_{q_j}+l_{q_j^{'}}\}\\
&\quad +\frac{1+\zeta-\eta}{1-\zeta-\eta}\ \max_{j=1,2,\ldots,N}\{|u_j-u_j^{'}|\},
\end{align}
where $\zeta=\max\limits_{1\leq j\leq N}\max\limits_{t\in I}\left\{|\alpha_j(t)|,|\beta_j(t)|\right\}$, $\eta=\max\limits_{1\leq j\leq N}\max\limits_{t\in I}\left\{|\gamma_j(t)|,|\delta_j(t)|\right\}$ and $\zeta+\eta<1.$
\end{theorem}

\begin{theorem}
Let $f=(f_1,f_2)$ and $f^{v}=(f_1^{v},f_2^{v})$ be the vector-valued FIF corresponding to the generalised data sets $V^{*}=\{(t_j,u_j,v_j):j=0,1,\ldots, N\}$ and $V_3^{*}=\{(t_j,u_j,v_j^{'}):j=0,1,\ldots, N\}$, respectively. Then,
\begin{align}
\nonumber\|f_1-f_1^{v}\|_{\infty}&<\frac{\zeta}{1-\zeta-\eta}\max_{j=0,1,\ldots,N}\{(l_{q_j}+l_{q_j^{v}})\}+\frac{1+\zeta-\eta+\zeta\eta}{1-\zeta-\eta}\max_{j=0,1,\ldots,N}\{|v_j-v_j^{'}|\}\\
&\quad +\frac{1-\eta}{1-\zeta-\eta}\max_{j=0,1,\ldots,N}\{(l_{p_j}+l_{p_j^{v}})\},
\end{align}
where $\zeta$ and $\eta$ are defined as earlier such that $\zeta+\eta<1.$
\end{theorem}

\subsection{Perturbation of the Vertical Scaling Factors}
\begin{theorem}
Let $f=(f_1,f_2)$ and $\tilde{f}=(\tilde{f_1},\tilde{f_2})$ be the vector-valued fractal functions interpolating the generalised data set $V^*$ corresponding to the iterated function systems $\mathcal{I}$ and $\tilde{\mathcal{I}}$, respectively. Assume $(1-\alpha)(1-\delta)-\beta\gamma~>~0$. Then,
\begin{align}\label{eq:vertscfperturb}
\nonumber\|f_1-\tilde{f_1}\|_\infty&<\frac{1-\delta}{(1-\alpha)(1-\delta)-\beta\gamma}\Big(\|s_j\|_{\infty} \Delta_{\alpha}+\|r_j\|_\infty\Delta_{\beta}+\Delta_{p}\Big)\\
&\quad+\frac{\beta}{(1-\alpha)(1-\delta)-\beta\gamma}\Big(\|s_j\|_{\infty}\Delta_\gamma+\|r_j\|_{\infty}\Delta_{\delta}+\Delta_{q}\Big).
\end{align}
\end{theorem}

\section{Stability of the Upper Box-Counting Dimension}\label{Dimension Results}

\subsection{Vertical scaling perturbation and upper box-counting dimension}
Consider the iterated function systems $\mathcal{I}$ and $\tilde{\mathcal{I}}$ with $f_1$ and $\tilde{f_1}$ being the corresponding EHVFIFs.
\begin{theorem}
Let $\tau=\max\limits_{1\leq j\leq N}\Big\{\frac{2\max\{\zeta_j,\eta_j\}}{|I_j|}\Big\}, \quad \tilde{\tau}=\max\limits_{1\leq j\leq N}\Big\{\frac{2\max\{\tilde{\zeta_j},\tilde{\eta_j}\}}{|I_j|}\Big\}$ and $\epsilon = \max\limits_{1\leq j\leq N} \Big\{\frac{2 \max\{\Delta_\alpha, \Delta_\beta, \Delta_\gamma, \Delta_\delta\}}{|I_j|}\Big\}$. Then, the following are true:
\begin{itemize}
\item If $\tau < 1$ and the perturbation is sufficiently small such that $\epsilon < 1 - \tau$, then the upper box-counting dimension is invariant under the perturbation, that is, $$\overline{\dim}_B(G(f_1)) = \overline{\dim}_B(G(\tilde{f_1})) = 1.$$
\item  If $\tau \geq 1$ and the perturbation is such that $\epsilon \leq \tau - 1 $, then $$1 \leq \overline{\dim}_B(G(f_1)) ,\quad \overline{\dim}_B(G(\tilde{f_1})) < 2.$$
\end{itemize}
\end{theorem}

\subsection{Data set Perturbation and Upper Box-Counting Dimension}
 \begin{theorem}
Let $\tau=\max\limits_{1\leq j\leq N}\Big\{\frac{2\max\{\zeta_j,\eta_j\}}{|I_j|}\Big\}$ and $\tau^{'}=\max\limits_{1\leq j\leq N}\Big\{\frac{2\max\{\zeta_j^{'},\eta_j^{'}\}}{|I_j^{'}|}\Big\}$. Let $f_1$ and $f_1^{t}$ be the EHVFIFs associated with the data sets $V^* = \{(t_j, u_j, v_j)\}$ and $V_1^* = \{(t'_j, u_j, v_j)\}$, respectively. Let $\epsilon=\max\limits_{j=1,\ldots,N}|\frac{2\max\{\zeta_j,\eta_j\}(|I_j^{'}|-|I_j|)}{|I_j||I_j^{'}|}|$. Then, the following are true:
\begin{itemize}
\item If $0<\tau<1$ and $\epsilon<1-\tau$, then
$$\overline{\dim}_B(G(f_1)) = \overline{\dim}_B(G(f_1^{t})) = 1.$$
\item If $\tau\geq1$ and $\epsilon<\tau-1$, then
$$1 \leq \overline{\dim}_B(G(f_1)) , \quad  \overline{\dim}_B(G(f_1^{t})) < 2.$$
\end{itemize}
\end{theorem}
\section*{Declarations}
\subsection*{Funding}
\noindent The first author received financial assistance from the Ministry of Education of India in the form of a Research Assistantship.
\subsection*{Data Availability}
\noindent No data is associated with this work.
\subsection*{Competing Interests}
\noindent The authors have no competing interests to declare.
\bibliographystyle{unsrt}
\bibliography{Stab}

@article{barnsley1986fractal,
	title={Fractal functions and interpolation},
	author={Barnsley, Michael F},
	journal={Constructive approximation},
	volume={2},
	pages={303--329},
	year={1986},
	publisher={Springer}
}

@article{barnsley1989hidden,
  title={Hidden variable fractal interpolation functions},
  author={Barnsley, M.F. and Elton, J. and Hardin, D. and Massopust, P.},
  journal={SIAM Journal on Mathematical Analysis},
  volume={20},
  number={5},
  pages={1218--1242},
  year={1989},
  publisher={SIAM}
}

@article{aiswarya2026,
title = {Construction and box-counting dimension of the Edelstein hidden variable fractal interpolation function},
journal = {Communications in Nonlinear Science and Numerical Simulation},
volume = {158},
pages = {109790},
year = {2026},
issn = {1007-5704},
doi = {https://doi.org/10.1016/j.cnsns.2026.109790},
url = {https://www.sciencedirect.com/science/article/pii/S1007570426001516},
author = {Aiswarya Thekkeettil and Srijanani {Anurag Prasad}}
}

@article{Edelstein1962,
author = {Edelstein, M.},
title = {On Fixed and Periodic Points Under Contractive Mappings},
journal = {Journal of the London Mathematical Society},
volume = {s1-37},
number = {1},
pages = {74-79},
doi = {https://doi.org/10.1112/jlms/s1-37.1.74},
url = {https://londmathsoc.onlinelibrary.wiley.com/doi/abs/10.1112/jlms/s1-37.1.74},
eprint = {https://londmathsoc.onlinelibrary.wiley.com/doi/pdf/10.1112/jlms/s1-37.1.74},
year = {1962}
}

@article{manousopoulos2011parameter,
  title={Parameter identification of 1D recurrent fractal interpolation functions with applications to imaging and signal processing},
  author={Manousopoulos, Polychronis and Drakopoulos, Vassileios and Theoharis, Theoharis},
  journal={Journal of Mathematical Imaging and Vision},
  volume={40},
  number={2},
  pages={162--170},
  year={2011},
  publisher={Springer}
}

@incollection{banerjee2020multifractal,
  title={Multifractal Analysis and Wavelet Decomposition in EEG Signal Classification},
  author={Banerjee, Santo and Easwaramoorthy, D and Gowrisankar, A},
  booktitle={Fractal Functions, Dimensions and Signal Analysis},
  pages={79--118},
  year={2020},
  publisher={Springer}
}

@article{raubitzek2021fractal,
  title={A fractal interpolation approach to improve neural network predictions for difficult time series data},
  author={Raubitzek, Sebastian and Neubauer, Thomas},
  journal={Expert Systems with Applications},
  volume={169},
  pages={114474},
  year={2021},
  publisher={Elsevier}
}

@article{chen2011reconstruction,
  title={The reconstruction of satellite images based on fractal interpolation},
  author={Chen, Ching-Ju and Cheng, Shu-Chen and Huang, YM},
  journal={Fractals},
  volume={19},
  number={03},
  pages={347--354},
  year={2011},
  publisher={World Scientific}
}
\end{document}